# VALIDITY OF THE EXPECTED EULER CHARACTERISTIC HEURISTIC

By Jonathan Taylor, Akimichi Takemura and Robert J. Adler

*Stanford University, University of Tokyo and Technion*

We study the accuracy of the expected Euler characteristic approximation to the distribution of the maximum of a smooth, centered, unit variance Gaussian process $f$. Using a point process representation of the error, valid for arbitrary smooth processes, we show that the error is in general exponentially smaller than any of the terms in the approximation. We also give a lower bound on this exponential rate of decay in terms of the maximal variance of a family of Gaussian processes $f^x$, derived from the original process $f$.

**1. Introduction.** In this paper, we study the expected Euler characteristic approximation to

$$P\left(\sup_{x \in M} f(x) \geq u\right) \tag{1}$$

where $f$ is the restriction to $M$ of $\widehat{f}$, a $C^2$ process on a $C^3$ manifold $\widehat{M}$, and $M$ is an embedded piecewise $C^3$ submanifold of $\widehat{M}$.

When the process $\widehat{f}$ is Gaussian with zero mean and has unit variance, the expected Euler characteristic approximation is given by

$$\widehat{P}\left(\sup_{x \in M} f(x) \geq u\right) = E(\chi(M \cap \widehat{f}^{-1}[u, +\infty)))$$
$$= \sum_{j=0}^{\dim M} \mathcal{L}_j(M)(2\pi)^{-(j+1)/2} \int_u^\infty H_j(r) e^{-r^2/2} \, dr, \tag{2}$$

where $H_j$ is the $j$th Hermite polynomial, $\chi(M \cap \widehat{f}^{-1}[u, +\infty))$ is the Euler characteristic of the excursion of $f$ above the level $u$ and the $\mathcal{L}_j(M)$ are the intrinsic volumes, or Lipschitz–Killing curvatures of the parameter space $M$ [20, 22], measured with respect to a Riemannian metric induced by $f$, which is discussed below in Section 2. In [20] only the special case of finite Karhunen–Loève processes (see below) was treated and in [22] the case of manifolds with smooth boundaries was treated. The result, proved in the generality described above, is to appear in [4]. With a slight abuse of notation, we have used $\widehat{\mathrm{P}}$ to denote our approximation to (1); that is, we are using $\widehat{\mathrm{P}}$ in the statistical sense of a "point estimator" and not as some alternative probability measure.

As noted above, the case when $\widehat{M}$ is $C^\infty$ and $\widehat{f}$ is a centered, unit variance, finite Karhunen–Loève expansion process was studied in [20]. These assumptions imply that there exist a map $\varphi\colon \widehat{M} \to S(\mathbb{R}^n)$ where $S(\mathbb{R}^n)$ is the unit sphere in $\mathbb{R}^n$ and a random vector $\xi(\omega) \sim N(0, I_{n\times n})$ such that

$$\widehat{f}(x,\omega) = \langle \varphi(x), \xi(\omega)\rangle_{\mathbb{R}^n} = \sum_{j=1}^n \xi_j(\omega)\varphi_j(x).$$

In this setting, without loss of generality, we can assume that $\widehat{M}$ is an embedded submanifold of $S(\mathbb{R}^n)$ and $\varphi$ is just the inclusion map. Using the volume of tubes approach [18, 23], it was shown in [20] that if $M$ is a piecewise $C^\infty$ submanifold of $S(\mathbb{R}^n)$ the error in the above approximation is bounded by

$$
\begin{aligned}
&\left| \mathrm{P}\left( \sup_{x \in M} f(x) \geq u \right) - \widehat{\mathrm{P}}\left( \sup_{x \in M} f(x) \geq u \right) \right| \\
&\quad \leq \frac{C}{\Gamma(n/2) 2^{(n-2)/2}} \int_{u/\cos\theta_c(M)}^\infty w^{n-1} e^{-w^2/2}\, dw \\
&\quad = C \times \mathrm{P}(\chi_n^2 \geq u^2/\cos^2\theta_c(M))
\end{aligned}
\tag{3}
$$

where $\theta_c(M)$ is a geometric quantity known as the critical radius of $M$ [11, 12, 18, 21]. For $M \subset S(\mathbb{R}^n)$, the critical radius is roughly defined by the following property: for every $0 < \theta < \theta_c(M)$ and arbitrary $z \in S(\mathbb{R}^n)$

$$d_{S(\mathbb{R}^n)}(z, M) \leq \theta \implies \arg\min_{x \in M} d_{S(\mathbb{R}^n)}(z, x) \text{ is unique}$$

where the metric on $S(\mathbb{R}^n)$ is the geodesic metric

$$d_{S(\mathbb{R}^n)}(x, z) = \cos^{-1}(\langle x, y\rangle).$$

In another setting, when $\widehat{f}$ is "almost" isotropic on $\mathbb{R}^k$, then, with some additional assumptions on $M$ (cf. Theorem 4.5.2 in [3]) Piterbarg [17] showed using the "double-sum" method that the error in using the expected Euler characteristic approximation is bounded by

$$\left| \mathrm{P}\left(\sup_{x\in M} f(x) \geq u\right) - \widehat{\mathrm{P}}\left(\sup_{x\in M} f(x) \geq u\right)\right| \leq C e^{-\alpha u^2/2} \tag{4}$$

for some $\alpha > 1$, though no expression for $\alpha$ is given. In the one-dimensional case, more was previously known [5, 13, 16] (see below).

Note that both bounds show that the error in approximating (1) is *exponentially* smaller than all terms in the expected Euler characteristic approximation. While undeniably useful, these two situations do not cover all possibilities. Referring to (4), certainly not every smooth Gaussian process of interest is isotropic, nor are the conditions required of $M$ easily interpretable (cf. Definition 4.5.1 in [3]). Referring to (3), while every Gaussian process does admit an *infinite* orthogonal expansion, see [2]

$$(5) \qquad f(x,\omega) = \sum_{j=1}^{\infty} \xi_j(\omega)\varphi_j(x)$$

through its reproducing kernel Hilbert space (RKHS), it is clear that substituting $n = \infty$ into (3) is meaningless. In fact, the situation is even worse in that the two cases do not even overlap: an isotropic field restricted to a bounded domain $T \subset \mathbb{R}^k$ *cannot* have a finite Karhunen–Loève expansion [19]!

This brings us to the main result of this work, Theorem 4.3, which, when $f$ is a constant variance Gaussian process as in the works cited above, provides bounds for the error in using the expected Euler characteristic approximation to (1). Specifically, when $f$ has unit variance, we show that

$$(6) \qquad \liminf_{u \to \infty} -u^{-2} \log \left| \mathrm{P}\left( \sup_{x \in M} f(x) \geq u \right) - \widehat{\mathrm{P}}\left( \sup_{x \in M} f(x) \geq u \right) \right| \geq \frac{1}{2} + \frac{1}{2\sigma_c^2(f)}.$$

Above, the "critical variance" $\sigma_c^2(f)$ depends on the variance of an auxiliary family of Gaussian processes $(f^x)_{x \in M}$, defined in (27) below.

An alternative approximation to (1) is to use the expected number of (extended outward) local maxima [10]. The term "approximation" is used in a somewhat loose sense, as, to the authors' knowledge, there are no generally applicable known closed form expressions for the number of extended outward local maxima of a smooth process. The only results known are ones which relate the asymptotic behavior of the expected number of local maxima of a Gaussian field on a manifold *without* boundary (which renders the qualifier "extended outward" unnecessary) to the expected Euler characteristic approximation, see [10]. Nevertheless, *if* one could compute the expected number of local maxima *exactly*, as one can the expected Euler characteristic in certain cases, one might expect to get a better approximation to (1). Virtually identical arguments to those used in this paper show that when $M$ is a manifold without boundary and $f$ is Gaussian with constant variance, on an exponential scale the errors in the approximations are equivalent, though, in the interest of brevity, we do not pursue this here. When the manifold $M$ has a boundary, the situation is more subtle and it may indeed be the case that the expected number of extended outward local maxima may be more accurate on an exponential scale.

The critical variance $\sigma_c^2(f)$ is closely related to the critical radius appearing in (3). Specifically, when $f$ is a centered, unit variance, finite Karhunen–Loève expansion process, and $M$ is a manifold without boundary, it is proven in Lemma 5.1 that

$$\sigma_c^2(f) = \cot^2 \theta_c(M)$$

where $\theta_c(M)$ is the critical radius of $M$, mentioned above and used in [21].

We note that, while (3) *is* an *explicit bound* for finite Karhunen–Loève expansion Gaussian processes, it is not sharp, nor generally applicable. In particular it depends on the (generally unknown) dimension of the sphere into which $M$ is embedded, that is, the dimension of the sphere in which $M$ sits. In a companion paper [14], when $f$ is a centered, unit variance, finite Karhunen–Loève expansion Gaussian process, the asymptotic error as $u \to \infty$ is evaluated using a Laplace approximation, rather than just the exponential behavior, which is the topic of this paper. In some one-dimensional stationary cases, the exact asymptotics of the error were found by Piterbarg in [16]. The results were generalized in [5], enlarging the class of processes covered by Piterbarg's result. Roughly speaking, the results of [5, 16] hold when the critical variance $\sigma_c^2(f)$ of the stationary process $f$ is achieved locally. In this case, the critical variance is explicitly computable in terms of the spectral moments of $f$ (see Lemma 5.3 below).

Our main result, Theorem 4.3, is formally an application of Theorem 3.3 to the case when $f$ is Gaussian with constant variance. Theorem 3.3 gives a bound for the error of the expected Euler characteristic approximation for the restriction of an *arbitrary* suitably regular (cf. [1, 22]) process $\widehat{f}$ on a $C^3$ manifold $\widehat{M}$ to any embedded piecewise $C^2$ submanifold $M \subset \widehat{M}$. Theorem 3.3 is, to the authors' knowledge, the only available bound for the error in the expected Euler characteristic approximation for arbitrary, suitably regular, smooth random fields and should prove useful in studying the accuracy of the Euler characteristic approximation to non-Gaussian fields [8, 9, 24]. The analogy to (3) for nonconstant variance Gaussian fields using a variant of the volume of tubes approach is presented in [15].

Another noteworthy feature to our approach is that it is a *direct* approach to determining the error in using the Euler characteristic approximation. This should be contrasted with the bounds (3) and (4) which were both arrived at indirectly in the sense that the bounds were derived for the "volume of tubes" approach and the "double-sum" approach and subsequently shown to hold for the expected Euler characteristic approximation.

The proof of Theorem 3.3 depends on a point set representation for the global maximizers of $h = \widehat{h}_{|M}$, the restriction to $M$ of a smooth deterministic function $\widehat{h}\colon \widehat{M} \to \mathbb{R}$, above the level $u$. Of course, there is a trivial, and not very useful, point set representation of the set of maximizers of $h$ above the level $u$:

$$\left\{x \in M : h(x) = \max_{y \in M} h(y), h(x) \geq u\right\}.$$

Lemma 2.2 gives an alternative point set representation of the maximizers using an auxiliary family of functions $(h^x)_{x \in M}$. Once we have a point set representation of the maximizers of a smooth function, we apply a "meta-theorem" for the density of point processes arising from smooth processes of [1, 3, 4, 7, 22]. If the global maximizer of the process $f$ is almost surely unique, the total mass of the density of the point process of maximizers above the level $u$ is therefore just (1). The auxiliary family of processes $(f^x)_{x \in M}$ mentioned above is, in some sense, the stochastic analogue of $(h^x)_{x \in M}$ in the deterministic setting.

After defining the processes $f^x$ and describing their properties, we derive the following almost sure bound:

$$|\mathbb{1}_{\{\sup_{x \in M} f(x) \geq u\}} - \chi(M \cap f^{-1}[u, +\infty))|$$
$$\leq \#\bigg\{x \in M : x \text{ is an extended outward critical point of } f,$$
$$f(x) \geq u, \sup_{y \in M \setminus \{x\}} f^x(y) > f(x)\bigg\}$$

where extended outward critical points are defined in a suitable fashion (cf. Section 2) and, for each $x \in M$, the process $f^x$ is uncorrelated with $f(x)$. Therefore, the points that contribute to the error are points where a Gaussian random variable $f(x)$ and the supremum of a process $f^x(y)$, independent of $f(x)$, are above the level $u$. The variance of the processes $f^x$ is what establishes the exponentially small relative error in (6). The only point where the Gaussian assumption is used is in bounding the expected number of points above, and the argument used here can be expected to extend to non-Gaussian processes as well.

The organization of the paper follows. Sections 2 lays out the regularity conditions needed for Theorem 3.3, and reviews some notions of piecewise smooth manifolds. Theorem 3.3 is proved in Section 3, and Section 4 deals with the unit variance Gaussian case, where Theorem 4.3 is proven. We conclude in Section 5 with some examples; specifically we compute $\sigma_c^2(f)$ for stationary processes on $\mathbb{R}$ and isotropic processes on $\mathbb{R}^k$ restricted to compact convex subsets.

**2. Suitably regular processes on piecewise $C^2$ manifolds.** In this section we describe the class of processes to which Theorem 3.3 will apply. Before setting out our assumptions, we recall some basic facts about piecewise $C^l$ ($l \geq 2$) submanifolds of an ambient $C^j$ ($j \geq l$) manifold $\widehat{M}$. For a more detailed self-contained treatment of the material in this section, specifically as it relates to smooth processes, the reader is referred to [4].

A $k$-dimensional piecewise $C^2$ manifold has a (not necessarily unique) decomposition into $C^2$ $i$-dimensional manifolds without boundary, $0 \leq i \leq k$; namely,

$$M = \bigcup_{i=0}^{k} \partial M_i.$$

Associated to every point $x \in M$ is its support cone in $T_x \widehat{M}$, the tangent space to $\widehat{M}$ at $x$

$$\mathcal{S}_x M = \{X_x \in T_x\widehat{M} : \exists \delta > 0, c \in C^1((-\delta, \delta), \widehat{M}),$$
$$c(t) \in M \ \forall t \in [0, \delta), c(0) = x, \dot{c}(0) = X_x\}.$$

In words, the support cone of $M$ at $x$ is the set of all directions in which a smooth curve can leave $x$ into $M$, but, in an infinitesimally small time period, still remain in $M$.

When $\widehat{M}$ is endowed with a Riemannian metric it is possible to define the dual cone of $\mathcal{S}_x M$. The following Riemannian metric will be essential to our analysis: an $L^2$ differentiable process induces a natural Riemannian metric on $\widehat{M}$ given by (cf. [22])

$$\widehat{g}_x(X_x, Y_x) \triangleq \mathrm{Cov}(X_x \widehat{f}, Y_x \widehat{f}),$$

for all $X_x, Y_x \in T_x\widehat{M}$.

The dual cone of $\mathcal{S}_x M$, $\mathcal{S}_x M^*$, in this case called the *normal cone* in $\widehat{M}$ at $x \in M$, is defined by

$$N_x M = \mathcal{S}_x M^* = \{X_x \in T_x\widehat{M} : \widehat{g}_x(X_x, Y_x) \leq 0, \ \forall Y_x \in \mathcal{S}_x M\}.$$

For $x \in \partial M_k = M^\circ$, $N_x M = (T_x M^\circ)^\perp$, the orthogonal complement of $T_x \partial M_k$ in $T_x\widehat{M}$. The normal cone figures prominently in the approximation, and in the main result of this paper, as both the Euler characteristic and the global maximizer point processes are defined in terms of *extended outward* critical points, that is, critical points at which the gradient (viewed as a tangent vector in the ambient space) is in the normal cone of the set at the critical point. Roughly speaking, this means that the function is increasing along curves leaving the set along certain "normal directions."

To avoid trivialities, we further assume our decomposition of $M$ is such that for every $0 \leq i \leq k$ and each $x \in \partial M_i$, the tangent space $T_x \partial M_i$ is the largest subspace contained in $\mathcal{S}_x M$. This condition rules out trivial decompositions of a nice open set $O \subset \mathbb{R}^2$ into $\{O \setminus F, F\}$ for a finite point set or some smooth closed curve $F \subset O$. This condition is not strictly necessary for the Morse theorem of [20]. In its place, we could impose a condition on the Morse functions to not have any extended outward critical points on such $F$'s. As this will almost surely be the case for the paths of suitably regular processes this point is somewhat moot.

For our purposes here, piecewise smooth manifolds are required to have the additional property that they are locally approximated by $\mathcal{S}_x M$ in the sense that, for every $x \in M$, there exists a diffeomorphism

$$(7) \qquad \varphi_x : U_x \subset \widehat{M} \to \varphi_x(U_x) \subset T_x \widehat{M}$$

of some neighborhood $U_x$ of $x \in \widehat{M}$, whose inverse, when restricted to $\mathcal{S}_x M$, maps any sufficiently small neighborhood of the origin to a neighborhood of $x \in M$. This condition, for instance, rules out cusps in $M$ as can occur when two manifolds intersect nontransversally. The importance of this condition, besides the fact that it is necessary in order to use the Morse theorem of [20], is that it implies that the set

$$(8) \qquad \{X_x \in N(\partial M_i) : X_x \in (P^{\perp}_{T_x \partial M_i} N_x M)^{\circ}\}$$

is open in $N(\partial M_i)$ where $(P^{\perp}_{T_x \partial M_i} N_x M)^{\circ}$ is the (relative) interior of $P^{\perp}_{T_x \partial M_i} N_x M$ in $T_x \partial M_i^{\perp}$ and $N(\partial M_i)$ is the normal bundle of $\partial M_i$ in $\widehat{M}$.

Having defined piecewise smooth manifolds, we now set out the assumptions on the processes $f$ defined on our piecewise smooth manifolds.

ASSUMPTION 2.1. We assume that $f$ is the restriction of $\widehat{f}$ to $M$, where $\widehat{f}$ is a square-integrable $C^2$ process on $\widehat{M}$, a $C^3$ $q$-dimensional manifold, and $M$ is a compact, embedded piecewise $C^2$ $k$-dimensional submanifold of $\widehat{M}$ such that $\mathcal{S}_x M$ is a convex cone for each $x \in M$. We further assume that, for each $i$, the gradient of $f_{|\partial M_i}$ read off in some nonrandom orthonormal frame field $E_i = (X_1, \ldots, X_i)$ on $\partial M_i$ satisfies the conditions of Lemma 2.5 of [22]. We denote this $\mathbb{R}^i$-valued process by $\nabla f_{|\partial M_i, E_i}$. These conditions are satisfied if $\widehat{f}_{|\partial M_i}$ is suitably regular in the sense of [22]. Finally, we assume that

$$\rho(x, y) \triangleq \mathrm{Cor}(\widehat{f}(x), \widehat{f}(y)) = 1 \quad \iff \quad x = y.$$

At this point, it is probably worth describing the significance of the above assumptions. In point form, the significance of these assumptions are:

1. The requirements that $f = \widehat{f}_{|M}$ for a $C^2$ process on $\widehat{M}$ and that the support cone $\mathcal{S}_x M$ is convex for each $x \in M$ are necessary for the point process representation of the Euler characteristic and to apply the Morse theorem of [20].

2. The conditions on the process $\nabla f_{|\partial M_i, E_i}$ ensure that the expected number of critical points of $f_{|\partial M_i}$ is finite and the density of point processes based on the critical points has an explicit integral representation in terms of a point process "meta-theorem" [3, 22].
3. The condition on $\rho(x,y)$ ensures that the map $x \mapsto \widehat{f}(x)$ is an embedding of $\widehat{M}$ into $L^2(\Omega, \mathcal{F}, P)$ which rules out "global" singularities in the process.

2.1. *A point set representation for the global maximizers of deterministic functions.* Although our primary interest is with stochastic processes, we begin here with the deterministic case, stating and proving Lemma 2.2 which describes a point process representation of the global maximizers of a deterministic function $h$. However, along each stratum $\partial M_i$ the point process depends on a $C^2$ function defined on

$$\partial M_i \times M^- \stackrel{\Delta}{=} \partial M_i \times M \setminus \{(x,y) \in \partial M_i \times M : x = y\}$$

which is singular near the "diagonal" $\{(x,y) \in \partial M_i \times M : x = y\}$. We resolve this singularity in Corollary 2.6, and arrive at the point process representation for the global maximizers in Corollary 2.7.

LEMMA 2.2. *Suppose $h = \widehat{h}_{|M}$, the restriction of $\widehat{h} \in C^2(\widehat{M})$. Fix $x \in \partial M_i$, and choose $\alpha^x \in C^2(M, (-\infty, 1])$ such that*

$$\alpha^x(y) = 1 \implies h(x) = h(y), \qquad \alpha^x(x) = 1.$$

*Then, $x$ is a maximizer of $h$ above the level $u$ if, and only if, the following three conditions hold:*

(i) $h(x) \geq u$.
(ii) $\nabla \widehat{h}(x) \in N_x M$. *That is, $x$ is an extended outward critical point of $h$.*
(iii) $h(x) \geq \sup_{y \in M \setminus \{x\}} h^x(y)$, *where, for all $y \in M$*

$$h^x(y) \stackrel{\Delta}{=} \begin{cases} \dfrac{h(y) - \alpha^x(y) h(x)}{1 - \alpha^x(y)}, & \text{if } \alpha^x(y) \neq 1, \\ h(y), & \text{if } \alpha^x(y) = 1. \end{cases}$$

*Further, if $\nabla^2 \alpha^x(x)$ is nondegenerate, and $x$ is a critical point of $h_{|\partial M_i}$, then, for any $C^2$ curve $c : (-\delta, \delta) \to \partial M_i$ with $c(0) = x, \dot{c}(0) = X_x \in T_x \partial M_i$,*

$$(9) \quad \lim_{t \to 0} h^x(c(t)) = \frac{\nabla^2 h_{|\partial M_i}(x)(X_x, X_x) - \nabla^2 \alpha^x_{|\partial M_i}(x)(X_x, X_x) h(x)}{-\nabla^2 \alpha^x_{|\partial M_i}(x)(X_x, X_x)}.$$

REMARK. The condition that $\alpha^x(x) = 1$ ensures that for each $x \in M$, $x$ is a critical point of $\alpha^x$.

PROOF OF LEMMA 2.2. The condition $h(x) \geq u$ is self-evident. Suppose, then, that $x \in \partial M_i, 0 \leq i \leq k$, is a maximizer of $h$. Then $\nabla \widehat{h}(x) \in N_x M$; otherwise there exists a direction $X_x \in \mathcal{S}_x M$ such that $\widehat{g}_x(X_x, \nabla \widehat{h}(x)) > 0$ and $x$ cannot be a maximizer.

Because $x$ is a maximizer, for all $y$ such that $\alpha^x(y) < 1$ it follows that

$$\frac{h(y) - \alpha^x(y) h(x)}{1 - \alpha^x(y)} < h(x).$$

On the other hand, if $\alpha^x(y) = 1$, then, by choice of $\alpha^x$, $h(y) = h(x)$ which proves that

$$h(x) \geq \sup_{y \in M \setminus \{x\}} h^x(y).$$

To prove the reverse implication, assume that $x$ is an extended outward critical point of $h_{|\partial M_i}$ and

$$h(x) \geq \sup_{y \in M \setminus \{x\}} h^x(y).$$

Now suppose that $x$ is not a maximizer of $h$; then there exists $y \in M \setminus \{x\}$ such that

$$h(x) < h(y).$$

In particular, for such a $y$, our choice of $\alpha^x$ implies that $\alpha^x(y) < 1$. It follows that

$$h(x) < \frac{h(y) - \alpha^x(y) h(x)}{1 - \alpha^x(y)}$$

which is a contradiction.

The limit (9) follows from two applications of l'Hôpital's rule. Specifically, we note that $x$ is a critical point of $h_{|\partial M_i}$ by assumption and the properties of $\alpha^x$ imply that it must also be a nondegenerate critical point of $\alpha^x$. Therefore,

$$\lim_{t \to 0} \frac{h(c(t)) - \alpha^x(c(t)) h(x)}{1 - \alpha^x(c(t))} = \lim_{t \to 0} \frac{(d/dt)(h(c(t)) - \alpha^x(c(t)) h(x))}{d/dt\, (1 - \alpha^x(c(t)))}$$
$$= \lim_{t \to 0} \frac{(d^2/dt^2)(h(c(t)) - \alpha^x(c(t)) h(x))}{(d^2/dt^2)(1 - \alpha^x(c(t)))}.$$

Since $x$ is also a nondegenerate critical point of $1 - \alpha^x(\cdot)$ the conclusion now follows from the fact that, for any $\beta \in C^2(\partial M_i)$ and any $C^2$ curve $c: (-\delta, \delta) \to \partial M_i$ with $x = c(0)$ a critical point of $\beta \in C^2(\partial M_i)$ and $\dot{c}(0) = X_x$

$$\lim_{t \to 0} \frac{d^2}{dt^2} \beta(c(t)) = \nabla^2 \beta(x)(X_x, X_x). \qquad \square$$

Condition (iii) above will be crucial to our later results as it is the condition which determines whether a given critical point is indeed a maximizer of $h$. The condition is not quite "ready to use," as we will need to consider an analogy to the quantity

$$\widetilde{W}(x) \overset{\Delta}{=} \sup_{y \in M \setminus \{x\}} h^x(y) \tag{10}$$

as a function of $x$. However, it is not hard to show, again by two applications of l'Hôpital's rule, that

$$\nabla \widehat{h}(x) \in N_x M \iff \sup_{y \in M \setminus \{x\}} h^x(y) < \infty, \tag{11}$$

$$\nabla \widehat{h}(x) \in N_x M \implies \inf_{y \in M \setminus \{x\}} h^x(y) = -\infty. \tag{12}$$

In other words, the above quantity is only finite at extended outward critical points of $h$ and the process is singular even if $\nabla \widehat{h}(x) \in N_x M$.

Although neither Lemma 2.2 nor (11) is exactly what we will need later, they contain a somewhat simplified version of our later arguments.

2.2. *Continuity of $\widetilde{W}(x)$: the "blow-up" of $h$.* As previously noted, the function

$$(x, y) \mapsto h^x(y) \tag{13}$$

is singular when $y$ is near $x$ even though it is $C^2$ on $\partial M_i \times M^-$. In particular

$$\lim_{y \to x} h^x(y)$$

is undefined. Ultimately, for the point process representation, we are interested in the continuity of $\widetilde{W}(x)$ on $\partial M_i$, and this singularity can make the arguments a little delicate.

In this section our main result shows that $\widetilde{W}(x)$ is continuous. Our strategy is to "blow up" $M$ around a neighborhood of $x \in \partial M_i$, and to relate $h^x(y)$, for $x$ in this neighborhood, to a continuous function on the "blow-up." We use the term "blow-up" as we relate the function $h^x(y)$ to a "desingularized" version of $h$ on the blow-up of $M$ near $x$.

By assumption, the parameter set $M$ is locally approximated by its support cone and we can, without loss of generality, describe the blow-up of $M$ under the assumption that $M \subset \widehat{M} = \mathbb{R}^q$ and for some neighborhood $U_x$

$$U_x \cap M = \{x\} \oplus \mathcal{S}_x M \cap U_x.$$

We can assume this because we have assumed in (7) that every $x \in M$ has such a neighborhood, and to establish continuity of $\tilde{W}(x)$ only local properties (in $x$) of $h^x(y)$ are needed.

For $x \in \partial M_i$, the support cone $\mathcal{S}_x M$ contains the $i$-dimensional tangent space $T_x \partial M_i$. By our assumptions on the decomposition of $M$, this subspace is maximal in $\mathcal{S}_x M$. Therefore, we can decompose $\mathcal{S}_x M$ as

$$T_x \partial M_i \times K_x$$

where $K_x$ is a convex cone that contains no subspace.

For ease of exposition, for the moment we move to a simple Euclidean setting and replace $\mathcal{S}_x M$ by

$$\widehat{K} = L \times K \subset \mathbb{R}^q$$

with $L$ a subspace, and $K$ a convex cone containing no subspaces. We will describe the "blow-up" construction first in this scenario, and then move back to piecewise smooth spaces. In what follows, $h = \widehat{h}_{|\widehat{K}}$ will be the restriction of a generic function on $\mathbb{R}^q$ to $\widehat{\mathbb{R}^q}$. We first define the "blow-up" $B(\widehat{K})$ of $\widehat{K}$ along $L$ as the disjoint union of the spaces

$$\mathcal{X} = L \times (\widehat{K} \setminus \{0\}),$$
$$\partial \mathcal{X} = L \times (\widehat{K} \cap S(\mathbb{R}^q)).$$

Above, $\mathcal{X}$ should be thought of as the image of

$$L \times \widehat{K} \setminus \{(x,y) \in L \times \widehat{K} : x = y\},$$

under the map

(14) $$(x,y) \stackrel{\Psi}{\mapsto} (x, y - x)$$

and the second space $\partial \mathcal{X}$ as the "boundary" of $\mathcal{X}$. The boundary is attached as follows: a sequence of points $(x_n, y_n) \in \mathcal{X}$ converges to $(x_0, y_0) \in \partial \mathcal{X}$ if, and only if,

$$x_n \to x_0, \qquad \|y_n\| \to 0, \qquad \frac{y_n}{\|y_n\|} \to y_0.$$

This notion of convergence corresponds to a sequence $\Psi^{-1}(x_n, y_n) = (x_n, x_n + y_n)$ converging to the diagonal

$$\{(x,y) \in L \times \widehat{K} : x = y\}$$

along a well-defined direction.

REMARK. Identifying the tangent bundle $T(\mathbb{R}^q)$ with $\mathbb{R}^q \times \mathbb{R}^q$, we can think of $B(\widehat{K})$ as a subset of $T(\mathbb{R}^q)$, in which the $y$'s above are replaced with $X_x \in \mathcal{S}_x M \subset T_x \mathbb{R}^q$. When convenient, we will consider $B(\widehat{K})$ as either a subset of $\mathbb{R}^{2q}$ or a subset of $T(\mathbb{R}^q)$.

In Lemma 2.2, $(\alpha^x)_{x \in M}$ was a family of functions which did not necessarily arise as the partial map of a function $\alpha : \widehat{M} \times \widehat{M} \to (-\infty, 1]$. In some cases, particularly in the stochastic setting below, this is a natural assumption to make.

LEMMA 2.3. *Suppose $\alpha \in C^2(\mathbb{R}^q \times \mathbb{R}^q)$ is such that the partial map*
$$\alpha^x(y) = \alpha(x, y)$$
*satisfies the conditions of Lemma 2.2 at every $x \in L$ and such that the Hessian of the partial map $\alpha^x$ is nondegenerate at every $x \in L$. Then, any $\widehat{h} \in C^2(\mathbb{R}^q)$ maps to a continuous function $\widehat{h}^{\alpha, \widehat{K}}$ on $B(\widehat{K})$ as follows:*
$$\widehat{h}^{\alpha, \widehat{K}}(x, y) = \frac{\widehat{h}(x+y) - \alpha(x, x+y)\widehat{h}(x) - \langle \nabla \widehat{h}_x, y \rangle_{\mathbb{R}^q}}{1 - \alpha(x, x+y)}$$

*and for $(x, y) \in \partial \mathcal{X}$*
$$\lim_{(t,s) \to (x,y)} \widehat{h}^{\alpha, \widehat{K}}(t, s) = \widehat{h}(x) + \frac{\nabla^2 \widehat{h}(x)(y, y)}{-\nabla^2 \alpha(x)(y, y)}.$$

PROOF. Two applications of l'Hôpital's rule. □

The term
$$\frac{\langle \nabla \widehat{h}_x, y \rangle_{\mathbb{R}^q}}{1 - \alpha(x, x+y)}$$
above "resolves" the singularity along the diagonal in some sense. In effect, it forces every $x \in L$ to be a critical point of the map
$$\mathbb{R}^q \ni y \mapsto \frac{\widehat{h}(x+y) - \alpha(x, x+y) \cdot \widehat{h}(x) - \langle \nabla \widehat{h}(x), y \rangle_{\mathbb{R}^q}}{1 - \alpha(x, x+y)}.$$

Our motivation for introducing $\widehat{h}^{\alpha, \widehat{K}}$ is to describe the singularities in the function $h^x(y)$ at critical points $x$ of $h_{|L}$ (recall from above that $L$ takes the place of $T_x \partial M_i$ in our general cone $\widehat{K}$). We are therefore interested in points where $x$ is a critical point of $h_{|L} = \widehat{h}_{|L}$. Note that if $x$ is a critical point of $h_{|L}$, then, for all $y \in \mathbb{R}^q$

(15)
$$\begin{aligned} h^x(x+y) &= \frac{h(x+y) - \alpha(x, x+y) \cdot h(x) - \langle \nabla h(x), y \rangle_{\mathbb{R}^q}}{1 - \alpha(x, x+y)} \\ &\quad + \frac{\langle P_L^\perp \nabla h(x), y \rangle_{\mathbb{R}^q}}{1 - \alpha(x, x+y)} \\ &= \widehat{h}^{\alpha, \widehat{K}}(x, y) + \frac{\langle P_L^\perp \nabla h(x), y \rangle_{\mathbb{R}^q}}{1 - \alpha(x, x+y)} \end{aligned}$$

where $P_L^\perp$ represents orthogonal projection onto $L^\perp$, the orthogonal complement of $L$ in $\mathbb{R}^q$. The expression (15) indicates that at critical points $x$ of $h_{|L}$, $h^x(y)$ is the sum of a well-behaved term, $\widehat{h}^{\alpha,\widehat{K}}$ and a singular term. The above relation holds for all critical points of $h_{|L}$. However, for our later arguments we would like to have this relation hold for all $x \in L$ to avoid having to condition on a set of measure zero. We therefore redefine $h^x$ as follows:

$$(16) \quad h^x(y) \triangleq \begin{cases} \dfrac{h(y) - \alpha^x(y)h(x) - \langle P_L \nabla h(x), y\rangle_{\mathbb{R}^q}}{1 - \alpha^x(y)}, & \text{if } \alpha^x(y) \neq 1, \\ h(y), & \text{if } \alpha^x(y) = 1. \end{cases}$$

With this redefinition of $h$, relation (15) holds for all $x \in L$ and, for each critical point $x$ of $h_{|L}$ the two definitions of $h^x$ coincide. For the remainder of this section, we will use the definition (16).

LEMMA 2.4. *If $P_L^\perp \nabla h(x) \in K^*$, the dual of $K \subset L^\perp$, then for any bounded neighborhood $O_x \ni x$*

$$\sup_{y \in O_x \setminus \{x\}} h^x(y) < +\infty.$$

PROOF. If $P_L^\perp \nabla h(x) \in K^*$, then (15) implies that for all $y$

$$h^x(y) \leq \widehat{h}^{\alpha,\widehat{K}}(x, x+y)$$

and $\widehat{h}^{\alpha,\widehat{K}}$ is continuous in $y$, and therefore bounded on bounded sets. If, on the other hand, $P_L^\perp \nabla h(x) \notin K^*$, then there exists a unit vector $v \in \widehat{K}$ such that $\langle P_L^\perp \nabla h(x), v\rangle_{\mathbb{R}^q} > 0$ and $x + tv \in O_x$ for sufficiently small $t$. Relation (15) implies that

$$\lim_{t \downarrow 0} h^x(x + tv) = +\infty.$$

This follows from the fact that the numerator in the expression on the right-hand side of (15) is strictly positive of order $O(t)$ for $t$ small, while the denominator is of order $O(t^2)$. □

LEMMA 2.5. *If $P_L^\perp \nabla h(x) \in (K^*)^\circ$, then for any bounded neighborhood $O$ of the origin in $\mathbb{R}^q$*

$$\widetilde{W}_O(x) = \sup_{y \in \widehat{K} \cap (\{x\} \oplus O) \setminus \{x\}} h^x(y) = \sup_{v \in \widehat{K} \cap (O \setminus \{0\})} \left(\widehat{h}^{\alpha,\widehat{K}}(x,v) + \frac{\langle P_L^\perp \nabla h(x), v\rangle_{\mathbb{R}^q}}{1 - \alpha(x, x+v)}\right)$$

*is continuous at $x$.*

PROOF. We first note that by (15) the two suprema above are equal, and it suffices to consider the supremum on the right.

Consider a convergent sequence

$$(x, v_n(x))_{n \geq 0} \to (x, v^*(x))$$

in $B(\widehat{K})$ along which the supremum $\widetilde{W}_O(x)$ is approached. Then either $\|v_n(x)\|_{\mathbb{R}^q} > 0$ for all $n$ sufficiently large, or $\|v_n(x)\|_{\mathbb{R}^q} \to 0$. In the first case it is immediate that $(x, v^*(x))$ is in $\mathcal{X}$ and, as we will show in a moment, in the second case $(x, v^*(x)) \in \partial \mathcal{X} \cap L \times S(L)$, where $S(L)$ is the unit sphere in $L$. In other words, the limiting direction $v^*(x)$ is in $S(L)$. Further, if $\|v_n(x)\|_{\mathbb{R}^q} \to 0$, then it is easy to see that the sequence $(x, P_L v_n(x))$ also achieves the supremum $\widetilde{W}_O(x)$.

To see why $v^*(x)$ must be in $S(L)$ in the second case, suppose that $(x, v_n(x))$ converges to $(x, v^*(x))$ with $v^*(x) \in \widehat{K} \setminus L$. As $P_L^\perp \nabla h(x) \in (K^*)^\circ$ it follows that

$$\langle P_L^\perp \nabla h(x), v^*(x) \rangle_{\mathbb{R}^q} < 0.$$

Therefore, applying the same argument as in (15),

$$\lim_{n \to \infty} \frac{\langle P_L^\perp \nabla h(x), v_n(x) \rangle_{\mathbb{R}^q}}{1 - \alpha(x, x + v_n(x))} = -\infty,$$

which contradicts the assumption that

$$\lim_{n \to \infty} h^x(x + v_n(x)) = \widetilde{W}_O(x).$$

The above implies that for any $x \in L$, and any convergent sequence $(x, v_n(x))$ achieving $\widetilde{W}_O(x)$, we can either assume that $\|v_n(x)\|$ is bounded uniformly below, or that $v_n(x) \in L$ for all $n$. Continuity of $\widetilde{W}_O(x)$ now follows, as, for such sequences and $\varepsilon > 0$ sufficiently small

$$\sup_{y \in B(x, \delta(\varepsilon))} |\widehat{h}^{a, \widehat{K}}(x, v_n(x)) - B(\widehat{h}, \alpha, \widehat{K})(y, v_n(x))| < \varepsilon,$$

$$\sup_{y \in B(x, \delta(\varepsilon))} \left| \frac{\langle P_L^\perp \nabla h(x), v_n(x) \rangle_{\mathbb{R}^q}}{1 - \alpha(x, x + v_n(x))} - \frac{\langle P_L^\perp \nabla h(y), v_n(x) \rangle_{\mathbb{R}^q}}{1 - \alpha(y, y + v_n(x))} \right| < \varepsilon. \qquad \square$$

Returning to piecewise smooth spaces, the above arguments will need slight modifications. Specifically, the map $\Psi$, defined in (14), has no natural replacement candidate for a piecewise smooth space. However, as noted in a remark above, we can think of $B(\widehat{K})$ as a subset of the tangent bundle $T(\mathbb{R}^q)$ in which case $\Psi : \mathbb{R}^q \times \mathbb{R}^q \to T(\mathbb{R}^q)$. In the piecewise smooth setting, we must therefore replace the map $\Psi$ with a map $H : \widehat{M} \times \widehat{M} \to T(\widehat{M})$ such that for each $x \in \widehat{M}$,

(17) $$H(\{x\} \times \widehat{M}) \subset T_x \widehat{M}.$$

One of the key properties of $\Psi$ was that the sequence $\Psi(x_n, y_n)$ converges to a point in $\mathcal{X}$ as $n \to 0$ as long as $x_n \to x$ and the unit vector

$$\frac{y_n - x_n}{\|y_n - x_n\|_{\mathbb{R}^q}}$$

converges in $S(\mathbb{R}^q)$. We can replace this property of $\Psi$ by asking the following of $H$. For any $C^1$ curve $c:(-\delta, \delta) \to \widehat{M} \times \widehat{M}$ with $c(0) = (c_1(0), c_2(0)) = (x, x), \dot{c}_2(0) = X_x$, we require that

(18) $$\lim_{t \to 0} \frac{H(x, c_2(t)) - H(x, x)}{t} = X_x.$$

Given such an $H$, we can, as in (16), redefine $h^x$ for a $C^2$ function on a piecewise smooth space $M$. For $x \in \partial M_i$, we redefine $h^x$ as follows:

(19) $$h^x(y) \triangleq \begin{cases} \dfrac{h(y) - \alpha^x(y)h(x) - \widehat{g}(\nabla h_{|\partial M_i}(x), H(x, y) - \alpha^x(y)H(x, x))}{1 - \alpha^x(y)}, \\ \qquad \text{if } \alpha^x(y) \neq 1, \\ h(y), \qquad \text{if } \alpha^x(y) = 1. \end{cases}$$

With $h^x$ redefined, by working in suitably chosen charts, it is not difficult to prove the following.

COROLLARY 2.6. *Suppose that $\alpha \in C^2(\widehat{M} \times \widehat{M})$ is such that the partial map*

$$\alpha^x(y) = \alpha(x, y)$$

*satisfies the conditions of Lemma 2.2 at every $x \in M$ and such that the Hessian of the partial map $\alpha^x$ is nondegenerate at every $x \in M$. Further, suppose $H: \widehat{M} \times \widehat{M} \to T(M)$ satisfies (17), (18) and let $h^x$ be defined as in (19). Then, the first conclusion of Lemma 2.2 holds if we replace* (i) *with:*

(i') $h(x) - \widehat{g}(H(x, x), \nabla h_{|\partial M_i}(x)) \geq u$.

*The second conclusion of Lemma 2.2 reads: for every $x \in \partial M_i$ and any $C^2$ curve $c:(-\delta, \delta) \to \partial M_i$ with $c(0) = x, \dot{c}(0) = X_x$,*

(20) $$\lim_{t \to 0} h^x(c(t)) = \frac{\nabla^2(h(\cdot) - \widehat{g}(H^x(\cdot), \nabla h_{|\partial M_i}(x)))(X_x, X_x)}{-\nabla^2 \alpha^x_{|\partial M_i}(x)(X_x, X_x)}$$
$$- (h(x) - \widehat{g}(H^x(x), \nabla h_{|\partial M_i}(x))),$$

*where $H^x(\cdot)$ is the partial map $H^x(y) = H(x, y)$.*

*Furthermore, the function*

$$\widetilde{W}(x) \triangleq \sup_{y \in M \setminus \{x\}} h^x(y)$$

is, for each $0 \leq i \leq k$, continuous on the set

$$\{x \in \partial M_i : P^\perp_{T_x \partial M_i} \nabla h(x) \in (P^\perp_{T_x \partial M_i} N_x M)^\circ\},$$

where $P^\perp_{T_x \partial M_i}$ represents projection onto the orthogonal complement of $T_x \partial M_i$ in $T_x \widehat{M}$.

2.3. *A point process for the global maximizers of stochastic processes.* For the remainder of this work, we choose a fixed piecewise smooth space $M$ and a process $\widehat{f}$ on $\widehat{M} \supset M$ satisfying Assumption 2.1.

In this section, we describe process analogies of $h^x$, $\alpha$ and $H$ in the case for which $\widehat{h}$ is replaced with the smooth process $\widehat{f}$. Specifically, we take

$$\alpha(x,y) = \rho(x,y),$$

(21)
$$H(x,y) = F(x,y) \triangleq \sum_{j=1}^{q} \mathrm{Cov}(f(y), X_j \widehat{f}(x)) X_{j,x},$$

for some orthonormal frame field $(X_{1,x}, \ldots, X_{q,x})$ on $T_x \widehat{M}$.

The following is the stochastic analogy of Lemma 2.6, that is, the point process representation of the maximizers of $f$.

COROLLARY 2.7. *Under Assumption* 2.1, *almost surely, the maximizers of $f$ are isolated and the maximizers of $f$ are the points $x \in \partial M_i, 0 \leq i \leq k$, such that:*

(i) $\nabla f_{|\partial M_i}(x) = 0$;
(ii) $f(x) - \widehat{g}(F^x(x), \nabla f_{|\partial M_i}(x)) \geq u$;
(iii) $P^\perp_{T_x \partial M_i} \nabla \widehat{f}(x) \in (P^\perp_{T_x \partial M_i} N_x M)^\circ$;
(iv) $f(x) \geq \sup_{y \in M \setminus \{x\}} \widetilde{f}^x(y)$, *where*

(22) $$\widetilde{f}^x(y) \triangleq \begin{cases} \dfrac{f(y) - \rho(x,y)f(x) - \widehat{g}(F^x(y) - \rho(x,y)F^x(x), \nabla f_{|\partial M_i}(x))}{1 - \rho(x,y)}, \\ \qquad \text{if } \rho(x,y) \neq 1, \\ f(y), \quad \text{if } \rho(x,y) = 1. \end{cases}$$

NOTE. If the joint density of $\nabla \widehat{f}(x)$, read off in some orthonormal basis of $T_x \widehat{M}$, is bounded by some constant $K$ uniformly in $x \in \widehat{M}$, then, almost surely, there will be no critical points of $f_{|\partial M_i}$ such that $P^\perp_{T_x \partial M_i} \nabla \widehat{f}(x) \in \partial N_x M \subset T_x \partial M_i^\perp$. Therefore, almost surely, all global maximizers will be such that $P^\perp_{T_x \partial M_i} \nabla \widehat{f}(x)$ is in the relative interior of $N_x M$ in $T_x \partial M_i^\perp$. As for the proof of Lemma 2.10 below, the proof of this is reasonably standard, and follows along the lines of similar results in, for example, Chapters 11 and 12 of [4]. We therefore omit the details.

PROOF OF COROLLARY 2.7. The only part of the argument in Lemma 2.2 that needs to be modified is what happens when $\alpha^x(y) = \rho(x,y) = 1$. In the deterministic case, we assumed that $\alpha^x(y) = 1$ implied $h(x) = h(y)$. In the random case, we know that $\alpha^x(y) = 1$ implies $f(x) - \mathrm{E}(f(x)) = (f(y) - \mathrm{E}(f(y)))\sigma(x)/\sigma(y)$ almost surely, where

$$\sigma^2(x) = \mathrm{Var}(f(x)).$$

Almost surely, then, it is still true that if $x$ is a maximizer of $f$, then $f(x) \geq \widetilde{f}^x(y)$ for all $y$ such that $\rho^x(y) = 1$, since otherwise $x$ cannot be a maximizer. The reverse implication follows similarly. $\square$

With $\widetilde{f}$ defined as above, it is easy to see that Corollary 2.6 holds almost surely with $h^x$ replaced by $\widetilde{f}^x$, $\alpha$ by $\rho$ and $H$ by $F$.

2.4. *Point process representation for the difference between the expected EC and the true probability.* Our assumptions allow us to use the Morse theorem of [20] to express the expected Euler characteristic of the excursions $M \cap \widehat{f}^{-1}[u, +\infty)$ as integrals over $M$. The formula is not new, though we repeat it here for use in deriving bounds on the error in the Euler characteristic approximation.

What is new, and is crucial to the entire paper, is the exact expression in Proposition 2.9 for the supremum distribution (1).

PROPOSITION 2.8. *Under Assumption* 2.1, *and with the notation there,*

$$\widehat{\mathrm{P}}\left(\sup_{x \in M} f(x) \geq u\right)$$

(23)
$$= \mathrm{E}(\chi(M \cap \widehat{f}^{-1}[u, +\infty)))$$

$$= \sum_{i=0}^{k} \int_{\partial M_i} \mathrm{E}(\det(-\nabla^2 f_{|\partial M_i, E_i}(x)) \mathbb{1}_{A_x^{EC}} | \nabla f_{|\partial M_i, E_i}(x) = 0)$$

$$\times \varphi_{\nabla f_{|\partial M_i, E_i}(x)}(0) \, d\mathcal{H}_i(x)$$

*where* $\mathcal{H}_i$ *is an* $i$-*dimensional Hausdorff measure induced by* $\widehat{g}$, $\varphi_{\nabla f_{|\partial M_i, E_i}(x)}$ *is the density of* $\nabla f_{|\partial M_i, E_i}(x)$ *and*

$$A_x^{EC} = \{f(x) \geq u, \nabla \widehat{f}(x) \in N_x M\}.$$

Suitable regularity of the process $\nabla f_{|\partial M_i, E_i}(x)$ implies that the maximizers of $f$ are almost surely isolated, though it does not guarantee uniqueness. If $\widetilde{W}(x)$ were continuous when restricted to $\partial M_i$, Assumption 2.1 would allow us to apply the general point process Lemmas 2.4 and 2.5 of [22] to the point process representation of the maximizers in Corollary 2.7. The almost sure analogy of Corollary 2.6 shows that $\widetilde{W}(x)$ is not continuous, but it is continuous on the open set

$$\{x \in \partial M_i : P^\perp_{T_x \partial M_i} \nabla \widehat{f}(x) \in (N_x M)^\circ\}.$$

Further, we are only interested in its behavior on this set. Straightforward modifications of the above cited lemmas, which we omit, lead to the following representation for the supremum distribution.

PROPOSITION 2.9. *Suppose that, almost surely, $f$ has a unique maximum and that Assumption 2.1 holds. Furthermore suppose that, for every $x \in M$,*

$$P(\widetilde{W}(x) = f(x)) = P(\widetilde{W}(x) = u) = 0.$$

*Then,*

$$P\left(\sup_{x \in M} f(x) \geq u\right)$$

(24)
$$= \sum_{i=0}^{k} \int_{\partial M_i} E(\det(-\nabla^2 f_{|\partial M_i, E_i}(x)) \mathbb{1}_{A_x^{SUP}} | \nabla f_{|\partial M_i, E_i}(x) = 0)$$
$$\times \varphi_{\nabla f_{|\partial M_i, E_i}(x)}(0) \, d\mathcal{H}_i(x),$$

*where*

$$A_x^{SUP} = \{f(x) \geq u \vee \widetilde{W}(x), \nabla \widehat{f}(x) \in N_x M\}.$$

The discrepancy between the expected Euler characteristic approximation and the true supremum distribution is

$$\mathrm{Diff}_{f,M}(u) \stackrel{\Delta}{=} \widehat{P}\left(\sup_{x \in M} f(x) \geq u\right) - P\left(\sup_{x \in M} f(x) \geq u\right)$$

(25)
$$= \sum_{i=0}^{k} \int_{\partial M_i} E(\det(-\nabla^2 f_{|\partial M_i, E_i}(x)) \mathbb{1}_{A_x^{ERR}} | \nabla f_{|\partial M_i, E_i}(x) = 0)$$
$$\times \varphi_{\nabla f_{|\partial M_i, E_i}(x)}(0) \, d\mathcal{H}_i(x)$$

where

$$A_x^{ERR} = \{u \leq f(x) \leq \widetilde{W}(x), \nabla \widehat{f}(x) \in N_x M\}.$$

Before concluding this section, we provide a lemma giving sufficient conditions for the uniqueness of the global maximum of $f$. As mentioned above, the proof is reasonably standard fare and so omitted. Detailed arguments for very similar cases can be found in Chapters 11 and 12 of [4]. In the nonmanifold setting, these arguments are classical. (E.g., see Theorem 3.2.1 of [1] or [6] and references therein.)

LEMMA 2.10. *Suppose that for all pairs $1 \leq i, j \leq k$ and all pairs $\{(x,y) : x \in \partial M_i, y \in \partial M_j\}$ the random vector*

$$V(x,y) = (f(x) - f(y), \nabla f_{|\partial M_i, E_i}(x), \nabla f_{|\partial M_j, E_j}(y))$$

*has a density, bounded by some constant $K(x,y)$. Then,*

$$\mathrm{P}(\{\exists\, (x,y) : x \in \partial M_i, y \in \partial M_j, V(x,y) = 0 \in \mathbb{R}^{i+j+1}\}) = 0.$$

**3. Bounding the error.** Expression (25) is an explicit formula for the error in the expected Euler characteristic approximation. A similar explicit formula can be derived for the error of the approximation based on the expected number of local maxima above the level $u$, though we do not pursue this here. However, as described in Section 2, the process $\widetilde{f}^x$ is singular near $x$, and actually has infinite variance near $x$, which means that standard tools such as the Borell–Tsirelson inequality cannot be used to bound its supremum distribution.

To see that the process $\widetilde{f}^x$ can have infinite variance, assume that $f$ is the restriction of an isotropic field to $[0,1]^2$. Fix a point $(x,0)$. We will compute the variance of $\widetilde{f}^x$ along the curve

$$c(t) = (x,0) - t \cdot (1,0).$$

Straightforward calculations show that

$$F^x(c(t)) = 0, \qquad F^x(x) = 0,$$

so that

$$\widetilde{f}^x(c(t)) = \frac{f(y) - \rho(x,y) f(x)}{1 - \rho(x,y)}.$$

In this case, the variance of $\widetilde{f}^x(y)$ is easily seen to be

$$\mathrm{Var}(\widetilde{f}^x(y)) = \frac{1 + \rho(x,y)}{1 - \rho(x,y)}$$

and

$$\lim_{y \to x} \mathrm{Var}(\widetilde{f}^x(y)) = \lim_{y \to x} \frac{1 + \rho(x,y)}{1 - \rho(x,y)} = +\infty.$$

In general, if $x \in \partial M_i$, then, along a curve $c : (-\delta, 0] \to M$ with $\dot{c}(0) = -X_x \in S_x M \setminus T_x \partial M_i$,

$$\lim_{t \uparrow 0} \mathrm{Var}(f^x(c(t))) = +\infty.$$

Although this is somewhat worrying, in (25) we only care about large *positive* values of $\widetilde{f}^x$, and, further, we only care about the behavior of $\widetilde{f}^x$ on the set

(26) $\qquad \{(x, \omega) : P^\perp_{T_x \partial M_i} \nabla \widehat{f}(x)(\omega) \in (P^\perp_{T_x \partial M_i} N_x M)^\circ\}.$

We exploit these facts and introduce a process $f^x$ in this section which has, under some conditions, finite variance and dominates $\widetilde{f}^x$ on the set (26). It is this process whose variance appears in the exponential bound for the behavior of $\text{Diff}_{f,M}(u)$ in the Gaussian case.

Obviously, the process $f^x$, which we define below, does not dominate the *absolute* value of $\widetilde{f}^x$. Indeed, if this were true, the process $f^x$ would have infinite variance as well.

The process $f^x$ is defined as follows:

$$(27) \quad f^x(y) \triangleq \begin{cases} \dfrac{f(y) - \rho^x(y)f(x) - \widehat{g}(\widehat{F}^x(y) - P_{N_xM}\widehat{F}^x(y), \nabla \widehat{f}(x))}{1 - \rho^x(y)}, \\ \qquad\qquad\qquad\qquad\qquad\qquad\qquad \text{if } \rho^x(y) \neq 1, \\ f(y) - \widehat{g}(\widehat{F}^x(y) - P_{N_xM}\widehat{F}^x(y), \nabla \widehat{f}(x)), \quad \text{if } \rho^x(y) = 1, \end{cases}$$

where $P_{N_xM}: T_x\widehat{M} \to N_xM$ represents orthogonal projection onto $N_xM$ and $\widehat{F}: \widehat{M} \times \widehat{M} \to T(\widehat{M})$ is given by

$$(28) \qquad \widehat{F}(x,y) = \begin{cases} F(x,y) - \rho(x,y)F(x,x), & \text{if } \rho(x,y) \neq 1, \\ F(x,y), & \text{if } \rho(x,y) = 1, \end{cases}$$

where $F$ is defined in (21).

LEMMA 3.1. *On the set* (26) *of extended outward normal points, for every $y \in M$,*

$$f^x(y) \geq \widetilde{f}^x(y).$$

*If $x \in \partial M_k = M^\circ$, then equality holds above.*

PROOF. First, we note that

$$(1 - \rho^x(y)) \cdot (\widetilde{f}^x(y) - f^x(y)) = \widehat{g}(\widehat{F}^x(y) - P_{N_xM}\widehat{F}^x(y), P^\perp_{T_x\partial M_i}\nabla \widehat{f}(x)).$$

As $N_xM$ is a convex cone, it follows that for any $Y_x \in T_x\widehat{M}$

$$Y_x - P_{N_xM}Y_x \in N_xM^*$$

where $N_xM^*$ is the dual cone of $N_xM$, which is just the convex hull of $\mathcal{S}_xM$. By duality,

$$\widehat{g}(Y_x - P_{N_xM}Y_x, V_x) \leq 0$$

for every $V_x \in N_xM$. Consequently, on the set (26)

$$\widehat{g}(Y_x - P_{N_xM}Y_x, P^\perp_{T_x\partial M_i}\nabla \widehat{f}(x)) \leq 0$$

for every $Y_x \in T_x\widehat{M}$. As $\widehat{F}^x(y) \in T_x\widehat{M}$ for each $y$, the first claim holds.

EXPECTED EULER CHARACTERISTIC        21As for the second, if $x \in \partial M_k$, then $N_x M = T_x \partial M_k^\perp$ and $P_{N_x M} V_x = 0$ for all $V_x \in T_x \partial M_k$. Similarly,
$$\widehat{g}(V_x, P^\perp_{T_x \partial M_i} \nabla \widehat{f}(x)) = 0.$$
Therefore, on this set
$$\widehat{g}(\widehat{F}^x(y) - P_{N_x M}\widehat{F}^x(y), \nabla \widehat{f}(x)) = 0. \qquad \square$$

As far as the continuity (in $x$) of
$$(29) \qquad W(x) \stackrel{\Delta}{=} \sup_{y \in M \setminus \{x\}} f^x$$
is concerned, it is not difficult to show that, almost surely, Corollary 2.6 holds with $\widetilde{f}^x$ replaced by $f^x$, that is, that $W(x)$ is continuous on the set
$$\{x \in \partial M_i : P^\perp_{T_x \partial M_i} \nabla \widehat{f}(x) \in (P^\perp_{T_x \partial M_i} N_x M)^\circ\}.$$

As we will see in the proof of Theorem 3.3, Lemma 3.1 provides the basic bounds for $\mathrm{Diff}_{f,M}(u)$. The following corollary to Lemma 2.2 will also be of use to us.

COROLLARY 3.2. *If $f$ has unit variance, then, for any $C^1$ unit speed curve $c:(-\delta, \delta) \to \partial M_k$ with $c(0) = x, \dot{c}(0) = X_x$,*
$$\lim_{t \to 0} f^x(c(t)) = -\frac{\nabla^2 f_{|\partial M_k}(x)(X_x, X_x) - \nabla^2 \rho^x(x)(X_x, X_x) f(x)}{-\nabla^2 \rho^x(x)(X_x, X_x)}$$
$$= -\nabla^2 f_{|\partial M_k}(x)(X_x, X_x) + f(x).$$
*Further,*
$$(30) \qquad \sup_{X_x \in S(T_x \partial M_k)} |-\nabla^2 f_{|\partial M_k}(x)(X_x, X_x) + f(x)| \leq \sup_{y \in M \setminus \{x\}} |f^x(y)|.$$

PROOF. The proof is essentially just the second conclusion of Lemma 2.2, recast in the stochastic process framework. The only thing that needs to be verified is that
$$\nabla^2 \rho^x(x)(X_x, X_x) = -1,$$
but this follows from the fact that
$$\nabla^2 \rho^x(y)(X_y, Y_y) = \mathrm{Cov}(\nabla^2 f_{|\partial M_k}(y)(X_y, Y_y), f(x))$$
and the fact that, as a double form
$$\mathrm{Cov}(\nabla^2 f_{|\partial M_k}(x)(X_x, Y_x), f(x)) = -\widehat{g}_x(X_x, Y_x),$$
(cf. [22]). $\square$

Using the results of Lemma 3.1 we have the following theorem.

THEOREM 3.3. *Suppose that $f$ has a unique maximum, almost surely, and Assumption* 2.1 *holds. Further suppose that, for every $x \in M$,*

$$P(W(x) = f(x)) = P(W(x) = u) = 0.$$

*Then,*

$$|\text{Diff}_{f,M}(u)|$$

(31)
$$\leq \sum_{i=0}^{k} \int_{\partial M_i} E(|\det(-\nabla^2 f_{|\partial M_i, E_i}(x))| \mathbb{1}_{B_x^{ERR}} |\nabla f_{|\partial M_i, E_i}(x) = 0)$$

$$\times \varphi_{\nabla f_{|\partial M_i, E_i}(x)}(0) \, d\mathcal{H}_i(x)$$

*where*

$$B_x^{ERR} = \{u \leq f(x) \leq W(x)\}.$$

**4. Gaussian fields with constant variance.** In this section, using Theorem 3.3, we derive an explicit bound for the exponential behavior of $\text{Diff}_{f,M}(u)$ when $\widehat{f}$ is a Gaussian field with constant variance, satisfying Assumption 2.1. The assumption of constant variance implies certain random variables are uncorrelated, hence independent in the Gaussian case. In particular, the assumption of constant variance implies that for $x \in \partial M_i, 0 \leq i \leq k$, the entire process $(f^x(y))_{y \in M \setminus \{x\}}$ is independent of $f(x)$ as well as $\nabla f_{|\partial M_i}(x)$. This allows us to remove the conditioning on $\nabla f_{|\partial M_i}(x)$ below. Once this conditioning is removed, the rest of the argument relies only on the Borell–Tsirelson inequality [2].

Our first observation is that, whether $f$ has constant variance or not, for each $x \in M$, the process $f^x(y)$ is uncorrelated with the random vector $\nabla f_{|\partial M_i}(x)$. Hence, in the Gaussian case, $f^x(y)$ is independent of $\nabla f_{|\partial M_i}(x)$.

LEMMA 4.1. *For every $x \in \partial M_i, 0 \leq i \leq k$ and every $y \in M \setminus \{x\}$*

$$\text{Cov}(f^x(y), X_x f) = 0$$

*for every $X_x \in T_x \partial M_i$.*

PROOF. We first note that, if $\rho^x(y) \neq 1$, then

$$(1 - \rho^x(y)) \text{Cov}(f^x(y), X_x f)$$
$$= \text{Cov}(f(y) - \rho^x(y) f(x), X_x f)$$
$$\quad - \text{Cov}(\widehat{g}(\widehat{F}^x(y) - P_{N_x M} \widehat{F}^x(y), \nabla \widehat{f}(x)), X_x f)$$
$$= \text{Cov}(f(y) - \rho^x(y) f(x), X_x f) - \widehat{g}(\widehat{F}^x(y) - P_{N_x M} \widehat{F}^x(y), X_x).$$

If, on the other hand, $\rho^x(y) = 1$, then

$$\operatorname{Cov}(f^x(y), X_x f) = \operatorname{Cov}(f(y), X_x f)$$
$$- \operatorname{Cov}(\widehat{g}(F^x(y) - P_{N_x M} F^x(y), \nabla \widehat{f}(x)), X_x f)$$
$$= \operatorname{Cov}(f(y), X_x f) - \widehat{g}(F^x(y) - P_{N_x M} F^x(y), X_x).$$

The conclusion will therefore follow once we prove, for every $y \in M$,

$$\operatorname{Cov}(f(y), X_x f) = \widehat{g}(F^x(y) - P_{N_x M} F^x(y), X_x),$$
$$\operatorname{Cov}(f(y) - \rho^x(y) f(x), X_x f) = \widehat{g}(\widehat{F}^x(y) - P_{N_x M} \widehat{F}^x(y), X_x).$$

As the two arguments are similar, we just prove the first equality. The map $F^x$ can be decomposed as follows:

$$F^x(y) = P_{T_x \partial M_i} F^x(y) + P^\perp_{T_x \partial M_i} F^x(y),$$

where

$$P_{T_x \partial M_i} F^x(y) = \sum_{j=1}^{i} \operatorname{Cov}(f(y), X_j f(x)) X_{j,x},$$

$$P^\perp_{T_x \partial M_i} F^x(y) = \sum_{j=i+1}^{q} \operatorname{Cov}(f(y), X_j f(x)) X_{j,x},$$

and the orthonormal basis $(X_{1,x}, \ldots, X_{q,x})$ is chosen so that the set $(X_{1,x}, \ldots, X_{i,x})$ forms an orthonormal basis for $T_x \partial M_i$ and $(X_{i+1,x}, \ldots, X_{q,x})$ forms an orthonormal basis for $T_x \partial M_i^\perp$, the orthogonal complement of $T_x \partial M_i$ in $T_x \widehat{M}$.

Further, because $\widehat{g}(X_x, V_x) = 0$ for every $X_x \in T_x \partial M_i$ and $V_x \in N_x M$, it follows that

$$P_{N_x M} F^x(y) = P_{N_x M} P^\perp_{T_x \partial M_i} F^x(y)$$

and for every $X_x \in T_x \partial M_i$

$$\widehat{g}(F^x(y) - P_{N_x M} F^x(y), X_x) = \widehat{g}(P_{T_x \partial M_i} F^x(y), X_x)$$
$$+ \widehat{g}(P^\perp_{T_x \partial M_i} F^x(y) - P_{N_x M} P^\perp_{T_x \partial M_i} F^x(y), X_x)$$
$$= \widehat{g}(P_{T_x \partial M_i} F^x(y), X_x)$$
$$= \sum_{j=1}^{i} \operatorname{Cov}(f(y), X_j f(x)) \widehat{g}(X_{j,x}, X_x)$$
$$= \operatorname{Cov}(f(y), X_x f). \qquad \square$$

As noted above, the independence between $\nabla f_{|\partial M_i}(x)$ and the process $f^x$ allows us to remove the conditioning on $\nabla f_{|\partial M_i}(x)$ in the expression for $\operatorname{Diff}_{f,M}(u)$, whether $f$ has constant variance or not.

COROLLARY 4.2. *Suppose $f$ is a Gaussian process satisfying the conditions of Theorem 3.3. Then,*

$$\text{(32)} \quad |\text{Diff}_{f,M}(u)| \leq \sum_{i=0}^{k} \int_{\partial M_i} \text{E}(|\det(-\nabla^2 f_{|\partial M_i, E_i}(x))| \mathbb{1}_{C_x^{ERR}}) \\ \times \varphi_{\nabla f_{|\partial M_i, E_i}(x)}(0) \, d\mathcal{H}_i(x)$$

*where*

$$C_x^{ERR} = \{u \leq f(x) - \widehat{g}(P_{T_x \partial M_i} F^x(x), \nabla \widehat{f}(x)) \leq W(x)\}.$$

*If $f$ has constant variance, then $F^x(x) = 0 \in T_x M$ and*

$$C_x^{ERR} = \{u \leq f(x) \leq W(x)\}.$$

PROOF. The only thing that needs to be proven is that, in the event $C_x^{ERR}$ the condition $\{u \leq f(x) \leq W(x)\}$ can be replaced with

$$\{u \leq f(x) - \widehat{g}(P_{T_x \partial M_i} F^x(x), \nabla \widehat{f}(x)) \leq W(x)\}$$

and the conditioning can be removed.

The reason that the above replacement is justified is that, on the set $\{\nabla f_{|\partial M_i}(x) = 0\}$

$$f(x) = f(x) - \widehat{g}(P_{T_x \partial M_i} F^x(x), \nabla \widehat{f}(x)).$$

Further, $f(x) - \widehat{g}(P_{T_x \partial M_i} F^x(x), \nabla \widehat{f}(x))$ is independent of $\nabla f_{|\partial M_i}(x)$, and Lemma 4.1 implies that $W(x)$ is also independent of $\nabla f_{|\partial M_i}(x)$. Therefore, the conditioning on $\nabla f_{\partial M_i}(x)$ can be removed. □

We are now ready to prove the following theorem.

THEOREM 4.3. *Let $\widehat{f}$ be a Gaussian process with constant, unit variance, on $\widehat{M}$ and let $f = \widehat{f}_{|M}$ be such that $f$ satisfies Assumption 2.1. Then,*

$$\liminf_{u \to \infty} -u^{-2} \log |\text{Diff}_{f,M}(u)| \geq \frac{1}{2}\left(1 + \frac{1}{\sigma_c^2(f)}\right)$$

*where*

$$\text{(33)} \quad \sigma_c^2(f,x) \stackrel{\Delta}{=} \sup_{y \in M \setminus \{x\}} \text{Var}(f^x(y))$$

*and*

$$\text{(34)} \quad \sigma_c^2(f) \stackrel{\Delta}{=} \sup_{x \in M} \sigma_c^2(f,x).$$

Proof. We must find an upper bound for (32). Writing

$$\nabla^2 f_{|\partial M_i, E_i}(x)$$
$$= \nabla^2 f_{|\partial M_i, E_i}(x) - \mathrm{E}(\nabla^2 f_{|\partial M_i, E_i}(x)|f(x)) + \mathrm{E}(\nabla^2 f_{|\partial M_i, E_i}(x)|f(x))$$
$$= \nabla^2 f_{|\partial M_i, E_i}(x) + f(x)I - f(x)I$$

(cf. [22]), and applying Hölder's inequality to (32) yields, for any conjugate exponents $p, q$,

$$|\mathrm{Diff}_{f,M}(u)| \leq \sum_{i=0}^{k} \int_{\partial M_i} \sum_{j=0}^{i} \mathrm{E}(f(x)^j \mathbb{1}_{\{f(x) \geq u\}})$$
$$\times \mathrm{E}(|\mathrm{detr}_{i-j}(-\nabla^2 f_{|\partial M_i, E_i}(x) - f(x)I)|^p)^{1/p}$$
$$\times \mathrm{P}(W(x) \geq u)^{1/q} \, d\mathcal{H}_i(x)$$

where $\mathrm{detr}_k(A)$ is the $k$th "det-trace" of the square matrix $A$ which is defined to be the sum of the determinants of all $k \times k$ principal minors of $A$.

Define

$$\mu \stackrel{\Delta}{=} \sup_{x \in M} E\left(\sup_{y \in M \setminus x} f^x(y)\right).$$

For

$$u \geq \mu \mathrm{E}(f^x(y)) + \mu^+$$

the Borell–Tsirelson inequality implies that

$$\mathrm{P}(W(x) \geq u) \leq 2e^{-(u-\mu^+)^2/2\sigma_c^2(f,x)}.$$

Recalling that $f^x(x) = f(x)$, for such $u$, it also follows that

$$\mathrm{E}(f(x)^j \mathbb{1}_{\{f(x) \geq u\}}) \leq C_j u^{j-1} e^{-(u-\mu^+)^2/2}.$$

Putting these facts together, for any conjugate exponents $p, q$

$$|\mathrm{Diff}_{f,M}(u)|$$
$$\leq C_k u^{k-1} e^{-((u-\mu^+)^2/2)(1+1/q\sigma_c^2(f))}$$
$$\times \sum_{i=0}^{k} \int_{\partial M_i} \sum_{j=0}^{i} \mathrm{E}(|\mathrm{detr}_{i-j}(-\nabla^2 f_{|\partial M_i, E_i}(x) - f(x)I)|^p)^{1/p} \, d\mathcal{H}_i(x).$$

The result now follows after noting that we can choose $q$ close to 1, and $u(q)$ so that, for $u \geq u(q)$, the remaining terms are arbitrarily small logarithmically, compared to $u^2$. □

Theorem 4.3 provides a lower bound on the exponential decay of $\text{Diff}_{f,M}(u)$. We believe that the lower bound is generally tight when a maximizer of $\sigma_c^2(f)$ occurs in $\partial M_k$, in the sense that the term corresponding to $\partial M_k$ in the sum defining $\text{Diff}_{f,M}(u)$ in (25) is exponentially of the same order as the upper bound; however, we were unable to prove this conjecture as it seems difficult to establish the sign of the error of the lower-order terms. In the piecewise smooth setting, it is therefore still open as to whether the $\liminf_{u \to \infty}$ in Theorem 4.3 can be replaced with $\lim_{u \to \infty}$ as we cannot rule out the possibility that some terms in the sum (25) cancel each other out, leading to a faster rate of exponential decay. Although we have not settled the issue completely, these situations seem somewhat pathological.

**5. Examples.** In this section, we compute $\sigma_c^2(f)$ for some simple examples, strengthening earlier results of [13, 16, 17]. Before turning to the examples, however, we discuss the relation between $\sigma_c^2(f)$ and the critical radius of a tube around $M$ when $f$ is assumed to be centered with unit variance. Specifically, we describe the geometry of the situation in the case of "global overlap," that is, when the supremum

$$\sigma_c^2(f) = \sup_{x \in M} \sup_{y \in M \setminus \{x\}} \text{Var}(f^x(y))$$

is achieved at a pair $(x^*, y^*)$, $x^* \neq y^*$.

5.1. *Geometric picture in the case of global overlap.* Here, we describe the notion of "global overlap" and describe the geometry of the process $f$ near pairs $(x^*, y^*)$ achieving the critical variance $\sigma_c^2(f)$. Roughly speaking, this situation occurs when $M$, the parameter space of $f$, "wraps around itself" and, for some $x \in M$ there is a point $y \in M$ that is close to $x$ in the $L^2$-metric but far in terms of geodesic distance from $x$. To describe the geometry involved in this situation we turn to spherical geometry in $\widetilde{H}_f$, the RKHS of $f$. Recall that $\widetilde{H}_f$ is defined by the reproducing kernel condition

$$\langle R(s, \cdot), R(t, \cdot) \rangle_{\widetilde{H}_f} = R(t, s)$$

and there exists an isometry that maps $\widetilde{H}_f$ onto the linear span $H_f$ of $\{f(x), x \in M\}$ in $L^2(\Omega, \mathcal{F}, P)$. Without loss of generality, then, we can describe the geometry in terms of $H_f$. Let $\Psi : M \to H_f \subset L^2(\Omega, \mathcal{F}, P)$ denote the map

$$x \mapsto f(x).$$

If $S(H_f)$ is the unit sphere in $H_f$, and $f$ is centered with unit variance, then $\Psi(M) \subset S(H_f)$, and our standing assumptions, namely that $f$ is $C^2$ and $\rho(x,y) \neq 1$ if $x \neq y$, imply that $\Psi$ is a piecewise $C^2$ embedding. Further, it is not hard to see that

$$\Psi_*(X_x) = X_x f,$$

so that the tangent space $T_{f(x)}\Psi(M)$ is spanned by $(X_{1,x}f, \ldots, X_{k,x}f)$ for some basis $\{X_{1,x}, \ldots, X_{k,x}\}$ of $T_x M$.

We denote the orthogonal complement of $T_{f(x)}\Psi(M)$ in $T_{f(x)}H_f$ by $T_{f(x)}^{\perp}\Psi(M)$, and the orthogonal complement of $T_{f(x)}\Psi(M)$ in $T_{f(x)}S(H_f)$ by $N_{f(x)}\Psi(M)$.

Given a point $f(x)$ and a unit normal vector $v_{f(x)} \in N_{f(x)}\Psi(M)$ we denote the geodesic, in $S(H_f)$, originating at $f(x)$ in the direction $v_{f(x)}$ by $c_{f(x),v_{f(x)}}$. That is,

$$c_{f(x),v_{f(x)}}(t) = \cos t \cdot f(x) + \sin t \cdot v_{f(x)}, \qquad 0 \leq t < \pi.$$

As discussed in [21], up to a certain point along $c_{f(x),v_{f(x)}}$, the points on $c_{f(x),v_{f(x)}}$ metrically project uniquely to $x$. That is, for $t$ small enough, the unique point on $\Psi(M)$ closest to $c_{f(x),v_{f(x)}}(t)$ is $f(x)$. We denote the largest $t$ for which this is true by $\theta(f(x), v_{f(x)})$ and we call it the *local critical radius (angle) at $f(x)$ in the direction $v_{f(x)}$*. Taking the supremum over all directions $v_{f(x)} \in S(N_{f(x)}\Psi(M))$ we obtain $\theta(f(x))$, the *local critical radius at $f(x)$*,

(35) $$\theta_c(f(x)) = \inf_{v_{f(x)} \in S(N_{f(x)}\Psi(M))} \theta(f(x), v_{f(x)})$$

and the global critical radius

(36) $$\theta_c = \theta_c(\Psi(M)) = \inf_{x \in M} \theta_c(f(x)).$$

The relation between these critical angles and $\sigma_c^2(f,x)$ in (33) is given in the following lemma.

LEMMA 5.1. *Suppose $f = \widehat{f}_{|M}$ is the restriction of a centered, unit variance Gaussian process on $\widehat{M}$, to $M$, a piecewise $C^2$ $k$-dimensional submanifold of $\widehat{M}$. Suppose that $\Psi : \widehat{M} \to H_f$ is a $C^2$ embedding. Then, for all $x \in \partial M_k$,*

$$\sigma_c^2(f,x) = \cot^2(\theta_c(x)).$$

PROOF. For $x \in \partial M_k$, $N_x M = T_x M^{\perp}$ is a linear space, therefore

$$\widehat{F}^x(y) - P_{N_x M}\widehat{F}^x(y) = P_{T_x \partial M_j}\widehat{F}^x(y).$$

Furthermore, because $f$ has constant variance $F^x(x) = 0$. Putting these facts together shows that

$$f^x(y) = \frac{f(y) - \rho^x(y)f(x) - \widehat{g}(\widehat{F}^x(y), \nabla \widehat{f}(x))}{1 - \rho^x(y)}$$

$$= \frac{f(y) - \rho^x(y)f(x) - \sum_{i=1}^{k} \mathrm{Cov}(f(y), X_i f(x)) X_i f(x)}{1 - \rho^x(y)}$$

for some orthonormal frame field $(X_1, \ldots, X_k)$ on $M$.

Turning to the picture in terms of geodesics, fix $f(x)$ and $v_{f(x)}$. Suppose that for a certain $t$ the point $c_{f(x),v_{f(x)}}(t)$ does not metrically project to $f(x)$. This implies there is a point $f(y) \in \Psi(M)$ such that

$$d(c_{f(x),v_{f(x)}}(t), f(y)) = \cos^{-1}(\langle c_{f(x),v_{f(x)}}(t), f(y) \rangle_{H_f})$$
$$< d(c_{f(x),v_{f(x)}}(t), f(x))$$
$$= \cos^{-1}(\langle c_{f(x),v_{f(x)}}(t), f(x) \rangle_{H_f}).$$

Alternatively,

$$\langle c_{f(x),v_{f(x)}}(t), f(y) \rangle_{H_f} = \cos t \cdot \langle f(x), f(y) \rangle_{H_f} + \sin t \cdot \langle v_{f(x)}, f(y) \rangle_{H_f}$$
$$= \cos t \cdot \rho(x,y) + \sin t \cdot \langle v_{f(x)}, f(y) \rangle_{H_f}$$
$$> \langle c_{f(x),v_{f(x)}}(t), f(x) \rangle_{H_f}$$
$$= \cos t.$$

After a little rearranging, we see that this is true if and only if

$$\cot t < \frac{\langle v_{f(x)}, f(y) \rangle_{H_f}}{1 - \rho(x,y)}.$$

Therefore,

$$\cot \theta_c(f(x), v_{f(x)}) = \sup_{y \in M \setminus \{x\}} \frac{\langle v_{f(x)}, f(y) \rangle_{H_f}}{1 - \rho(x,y)}.$$

Taking the supremum over all $v_{f(x)} \in N_{f(x)}\Psi(M)$ we see that

$$\cot \theta_c(f(x)) = \sup_{v_{f(x)} \in S(N_{f(x)}\Psi(M))} \sup_{y \in M \setminus \{x\}} \frac{\langle v_{f(x)}, f(y) \rangle_{H_f}}{1 - \rho(x,y)}$$
$$= \sup_{y \in M \setminus \{x\}} \frac{\|P_{N_{f(x)}\Psi(M)} f(y)\|_{H_f}}{1 - \rho(x,y)}$$

where $P_{N_{f(x)}\Psi(M)}$ represents orthogonal projection onto $N_{f(x)}\Psi(M)$. Therefore,

$$\cot^2\theta_c(f(x)) = \sup_{y\in M\setminus\{x\}} \frac{\|P_{N_{f(x)}\Psi(M)}f(y)\|^2_{H_f}}{(1-\rho(x,y))^2}$$

and it remains to show that

$$\|P_{N_{f(x)}\Psi(M)}f(y)\|^2_{H_f}$$
$$= \mathrm{Var}\left(f(y) - \rho(x,y)f(x) - \sum_{i=1}^k \mathrm{Cov}(f(y), X_i f(x))X_i f(x)\right).$$

This, however, follows from the fact that $N_{f(x)}\Psi(M)$ is the orthogonal complement (in $T_{f(x)}H_f$) of the subspace $L_x = \mathrm{span}\{f(x), X_1 f(x), \ldots, X_k f(x)\}$, and the fact that

$$\rho(x,y)f(x) + \sum_{i=1}^k \mathrm{Cov}(f(y), X_i f(x))X_i f(x)$$
$$= \langle f(x), f(y)\rangle_{H_f} f(x) + \sum_{i=1}^k \langle f(y), X_i f(x)\rangle_{H_f} X_i f(x)$$
$$= P_{L_x} f(y). \qquad \square$$

We now consider the case when $M$ is a manifold without boundary and the supremum in

(37) $$\cot^2(\theta_c) = \sup_{x\in M}\sup_{y\in M\setminus\{x\}} \frac{\|P_{N_{f(x)}\Psi(M)}f(y)\|^2_{H_f}}{(1-\rho(x,y))^2}$$

is achieved at some $(x^*, y^*)$. Thinking of the critical variance as the cotangent of some critical distance on $M$, then, geometrically, the tube of radius $\theta_c$ should self-intersect along a geodesic from $x^*$ to $y^*$ in such a way that the tube viewed locally from the point $x^*$ shares a hyperplane with the tube viewed locally from the point $y^*$. Alternatively, at the point of self-intersection the outward pointing unit normal vectors should be pointing in opposite directions.

The simplest way of seeing this is to think of $M$ as just two points $\{p_1, p_2\}$ in $\mathbb{R}^2$. In this case the tube of radius $r$ around $M$ consists of the union of two discs of radius $r$. When $r = d(p_1, p_2)/2$ the two discs self-intersect at exactly one point, and the unit normal vectors are pointing in opposite directions.

We can make this statement precise in the following lemma.

LEMMA 5.2. *Suppose $(x^*, y^*)$ achieve the supremum in (37). Then,*
$$f(y^*) = \cos(2\theta_c) \cdot f(x^*) + \sin(2\theta_c) \cdot v^*_{f(x^*)},$$
*where*
$$v^*_{f(x^*)} = \frac{P_{N_{f(x^*)}\Psi(M)} f(y^*)}{\|P_{N_{f(x^*)}\Psi(M)} f(y^*)\|_{H_f}}$$
*is the direction of the geodesic between $f(x^*)$ and $f(y^*)$. For any $X_{x^*} \in T_{x^*}M$ and any $X_{y^*} \in T_{y^*}M$,*
$$\mathrm{Cov}(X_{x^*}f, f(y^*)) = \langle \Psi_*(X_{x^*}), f(y^*) \rangle_{H_f} = 0,$$
$$\mathrm{Cov}(X_{y^*}f, f(x^*)) = \langle \Psi_*(X_{y^*}), f(x^*) \rangle_{H_f} = 0.$$
*Further, the partial map*
$$\rho^{x^*}(y) \triangleq \rho(x^*, y)$$
*has a local maximum at $y^*$, and the partial map*
$$\rho^{y^*}(x) \triangleq \rho(x, y^*)$$
*has a local maximum at $x^*$.*

PROOF. If $(x^*, y^*)$ achieve the supremum in (37), then there exists a point $z$ equidistant from $f(x^*)$ and $f(y^*)$ and unit vectors $v_{f(x^*)} \in N_{f(x^*)}\Psi(M)$ and $w_{f(y^*)} \in N_{f(y^*)}\Psi(M)$ such that
$$z = \cos\theta_c \cdot f(x^*) + \sin\theta_c \cdot v_{f(x^*)} = \cos\theta_c \cdot f(y^*) + \sin\theta_c w_{f(y^*)}.$$

It is not a priori obvious that $w_{f(y^*)} \in N_{f(y^*)}\Psi(M)$. However, if this were not the case, there exists $t \in M$, $t \neq y$, such that $d(t, z) < \theta_c$. This contradicts the assumption that $(x^*, y^*)$ achieve the supremum in (37) and the critical radius is $\theta_c$. The proof that this is a contradiction is left to the reader, though it follows similar lines to the argument at the end of this proof.

To prove the first claim, therefore, it is enough to show that the unit tangent vectors $\dot{c}_{f(x^*), v_{f(x^*)}}(\theta_c)$ and $\dot{c}_{f(y^*), w_{f(y^*)}}(\theta_c)$ at $z$ satisfy
$$u = \dot{c}_{f(x^*), w_{f(x^*)}}(\theta_c) + \dot{c}_{f(y^*), w_{f(y^*)}}(\theta_c) = 0 \in T_z H_f.$$

Suppose then that $u \neq 0$. A simple calculation shows that
$$\langle u, f(x^*) \rangle_{H_f} = \langle u, f(y^*) \rangle_{H_f} = \frac{\cos(2\theta_c) - \langle f(x^*), f(y^*) \rangle_{H_f}}{\sin\theta_c}.$$

If $u \neq 0$, then $f(x^*), f(y^*)$ and $z$ are not on the same geodesic and the triangle inequality implies that
$$\cos(2\theta_c) < \langle f(x^*), f(y^*) \rangle_{H_f}$$

which implies

$$\langle u, f(x^*)\rangle_{H_f} = \langle u, f(y^*)\rangle_{H_f} < 0.$$

Consider the geodesic originating at $z$ in the direction $u^* = u/\|u\|_{H_f}$

$$c_{z,u^*}(t) = \cos t \cdot z + \sin t \cdot u^*.$$

Assuming $u^* \neq 0$,

$$\langle \dot{c}_{z,u^*}(0), f(y^*)\rangle_{H_f} = \langle \dot{c}_{z,u^*}(0), f(x^*)\rangle_{H_f} = \langle u^*, f(x^*)\rangle < 0.$$

This implies that for sufficiently small $|s|$, $s < 0$,

$$\langle f(x^*), c_{z,u^*}(s)\rangle_{H_f} = \langle f(y^*), c_{z,u^*}(s)\rangle_{H_f} > \cos\theta_c.$$

For such an $s$, there exist distinct points $\hat{x}(s)$ and $\hat{y}(s)$ such that

$$d(c_{z,u^*}(s), \hat{x}(s)) < \theta_c, \qquad d(c_{z,u^*}(s), \hat{y}(s)) < \theta_c$$

such that the geodesic connecting $\hat{x}(s)$ and $c_{z,u^*}(s)$ is normal to $M$ at $\hat{x}(s)$, and the geodesic connecting $\hat{y}(s)$ and $c_{z,u^*}(s)$ is normal to $M$ at $\hat{y}(s)$.

Without loss of generality, we assume that

$$d(c_{z,u^*}(s), \hat{x}(s)) \leq d(c_{z,u^*}(s), \hat{y}(s)).$$

In this case, the geodesic connecting $\hat{y}$ to $c_{z,u^*}(s)$ is no longer a minimizer of distance once it passes the point $c_{z,u^*}(s)$, which is of distance strictly less than $\theta_c$ from $\hat{y}(s)$. That is, for some point $\hat{z}(s)$ along the geodesic connecting $\hat{y}(s)$ and $c_{z,u^*}(s)$, beyond $c_{z,u^*}(s)$ but of distance strictly less than $\theta_c$ from $\hat{y}(s)$, there exist points in $M$ strictly closer to $\hat{z}(s)$ than $\hat{y}(s)$. We therefore have a contradiction, as the critical radius of $\Psi(M)$ is $\theta_c$.

To prove the second claim, we note that $f(y^*)$ is a linear combination of $f(x^*)$ and $v_{f(x^*)}$ which are perpendicular to every vector in $T_{f(x^*)}\Psi(M)$. Similarly, $f(x^*)$ is a linear combination of $f(y^*)$ and $w_{f(y^*)}$ which are perpendicular to every vector in $T_{f(x^*)}\Psi(M)$. The fact that the partial maps $\rho^{x^*}(y)$ and $\rho^{y^*}(x)$ have local maxima follows from the same contradiction argument used above. $\square$

5.2. *Centered stationary processes on $[0,T]$.* In this section $f(x), x \in [0,T]$, is assumed to be a centered $C^2$ stationary process with unit variance and covariance function $R$. As in [13] we change the time scale so that $\mathrm{Var}(\dot{f}(x)) = -\ddot{R}(0) = 1$.

Fix a point $x \in (0,T)$, in which case the process $f^x$ is given by

$$f^x(y) = \frac{f(y) - R(x-y)f(x) - \dot{R}(x-y)\dot{f}(x)}{1 - R(x-y)}$$

and the critical variance at $x$ is given by

$$\sigma_c^2(f,x) = \sup_{y\in[0,T]\setminus\{x\}} \frac{\text{Var}(f(y)\mid f(x),\dot{f}(x))}{(1-R(x-y))^2}$$

$$= \sup_{\substack{-x\leq t\leq T-x \\ t\neq 0}} \frac{1-R(t)^2-\dot{R}(t)^2}{(1-R(t))^2}.$$

The critical variance in the interior is

(38) $$\sigma_c^2(f,(0,T)) = \sup_{x\in(0,T)} \sigma_c^2(f,x) = \sup_{0<t<T} \frac{1-R(t)^2-\dot{R}(t)^2}{(1-R(t))^2}.$$

A local version of critical variance $\sigma_{c,\text{loc}}^2(x,(0,T))$ is obtained by letting $t\to 0$ in (38), that is,

$$\sigma_{c,\text{loc}}^2(x,(0,T)) = \lim_{t\to 0} \frac{1-R(t)^2-\dot{R}(t)^2}{(1-R(t))^2} = R^{(4)}(0) - 1$$

where the last conclusion follows from some simple calculus. An alternative interpretation of the quantity $\sigma_{c,\text{loc}}^2(x,(0,T))$, is the following:

$$\sigma_{c,\text{loc}}^2(x,(0,T)) = \text{Var}(\ddot{f}(x)|f(x)).$$

The local critical variance at the end points $\{0,T\}$ needs slightly more attention. We consider the point $x=0$ without loss of generality. The normal cone at $x=0$ is

$$N_0[0,T] = \bigcup_{c\leq 0} c\frac{d}{dx} \subset T_0\mathbb{R}.$$

For ease of notation, we just write $N_0$ for $N_0[0,T]$.

The projection onto $N_0$ is just

$$P_{N_0}\left(a\frac{d}{dx}\right) = \mathbb{1}_{\{a\leq 0\}}\left(a\frac{d}{dx}\right).$$

As

$$\text{Cov}(\dot{f}(0),f(y)) = -\dot{R}(y)$$

the process $f^0(y)$ is given by

$$f^0(y) = \frac{f(y) - R(y)f(0) - \mathbb{1}_{\{\dot{R}(y)\leq 0\}}\dot{R}(y)\dot{f}(0)}{(1-R(y))^2}.$$

The local critical variance at $x=0$ is given by

(39)
$$\sigma_c^2(f,0) = \sup_{0\leq t\leq T} \frac{1-R(t)^2-\dot{R}(t)^2+\max(\dot{R}(t),0)^2}{(1-R(t))^2}$$

$$\geq \sup_{0\leq t\leq T} \frac{1-R(t)^2-\dot{R}(t)^2}{(1-R(t))^2}.$$

The inequality above implies that $\sigma_c^2(f,0) \geq \sigma_c^2(f,y), 0 \leq y \leq T$. Therefore the critical variance $\sigma_c^2(f)$ is attained at the end points $t = 0, T$.

Note that this does not exclude the case that the critical radius is also attained in the interior $0 < t < T$. Under some circumstances, the critical variance is attained everywhere, as demonstrated in the following lemma.

LEMMA 5.3. *Suppose $f$ is a centered, unit variance $C^2$ stationary process on $\mathbb{R}$ such that $-\ddot{R}(0) = 1$ and $\dot{R}(t) \leq 0$ for all $t \geq 0$. Then*

$$\sigma_c^2(f) = \sigma_{c,\mathrm{loc}}^2(f,x) = \mathrm{Var}(\ddot{f}(x)|f(x)).$$

PROOF. As noted above, we just have to compute $\sigma_c^2(f,0)$. Because $\dot{R}(t) \leq 0$ we see that

$$\sigma_c^2(f,0) = \sup_{0 \leq t \leq T} \frac{1 - R(t)^2 - \dot{R}(t)^2}{(1-R(t))^2}.$$

Now suppose that this supremum is achieved. Lemma 5.2 implies that for any $t$ that achieves this supremum, $\dot{R}(t) = 0$, and $t$ is a local maximum of $R(t)$. Let

$$T^* = \{t > 0 | \dot{R}(t) = 0, R(t) = \cos(2\theta_c)\}.$$

As $T^*$ is a closed set there exists a minimum value of $T^*$:

$$t^* = \min T^* > 0.$$

Because $t^*$ is a local maximum of $R$, there exists $\epsilon > 0$ such that

$$R(t) \leq R(t^*) \qquad \forall t \in (t^* - \epsilon, t^*).$$

But $R(t)$ is assumed to be nonincreasing and we have

$$R(t) \geq R(t^*) \qquad \forall t < t^*.$$

Therefore

$$R(t) = R(t^*) \qquad \forall t \in (t^* - \epsilon, t^*).$$

This, however, contradicts the minimality of $t^*$. □

REMARK. Presumably, the same ideas as in [5] could be used to generalize Lemma 5.3 to a larger class of covariance functions, though we will leave that as an exercise for the reader.

5.3. *Isotropic fields on $\mathbb{R}^m$ with monotone covariance functions.* In this section, we compute the critical variance for centered, unit variance isotropic processes restricted to a compact, convex set $M$. In particular, we prove the following proposition, which shows that the exponential behavior of $\text{Diff}_{f,M}(u)$ for such processes is determined solely by the conditional variance of the second derivative, given the field.

PROPOSITION 5.4. *Let $\widehat{f}$ be an isotropic process on $\mathbb{R}^m$ that induces the standard metric on $\mathbb{R}^m$ satisfying Assumption 2.1 and let $f = \widehat{f}_{|M}$ be the restriction of $\widehat{f}$ to a compact, convex set $M$ with piecewise smooth boundary. If the covariance function*

$$R(\|x\|) = \text{Cov}(f(x+t), f(t))$$

*is monotone nonincreasing, then the critical variance is attained locally and is given by*

$$\sigma_c^2(f) = \text{Var}\bigg(\frac{\partial^2 f}{\partial t_1^2}(0)\Big|f(0)\bigg).$$

PROOF. Fix $t \in M^\circ$. Similar computations to those in Section 5.2 show that

$$\text{Var}(f^t(s)) = \frac{1 - R^2(\|s-t\|) - \dot{R}^2(\|s-t\|)}{(1 - R(\|s-t\|))^2}.$$

Because $R$ is assumed monotone nonincreasing, the arguments used in Lemma 5.3 imply that

$$\sup_{s \in M\setminus\{t\}} \text{Var}(f^t(s)) = \text{Var}\bigg(\frac{\partial^2 f}{\partial t_1^2}(0)\Big|f(0)\bigg).$$

We therefore turn to the boundary $\partial M$, assumed to be piecewise smooth.

Fix $t \in \partial M$ and $s \neq t$. Let $X_t^s$ be the unit vector in $T_t\mathbb{R}^m$ in the direction $s - t$. Since $M$ is convex, $X_t^s \in \mathcal{S}_t M$ and

$$X_t^s \parallel \widehat{F}^t(s) \in \mathcal{S}_t M,$$

that is, $X_t^s$ is parallel to $\widehat{F}_s^t$, so that $P_{N_xM}\widehat{F}^t(s) = 0$. Therefore,

$$\text{Var}(f^t(s)) = \frac{1 - R^2(\|s-t\|) - \dot{R}^2(\|s-t\|)}{(1 - R(\|s-t\|))^2}.$$

We see that we do not incur any additional penalty for the local critical radius at the boundary points if $M$ is convex. Again, the arguments of Lemma 5.3 imply that

$$\sigma_c^2(f,t) = \text{Var}\bigg(\frac{\partial^2 f}{\partial t_1^2}(0)\Big|f(0)\bigg).$$

The conclusion now follows. □

J. TAYLOR
DEPARTMENT OF STATISTICS
SEQUOIA HALL
STANFORD UNIVERSITY
STANFORD, CALIFORNIA 94305-4065
USA
E-MAIL: jtaylor@stat.stanford.edu

A. TAKEMURA
GRADUATE SCHOOL OF INFORMATION SCIENCE
  AND TECHNOLOGY
UNIVERSITY OF TOKYO
7-3-1 HONGO, BUNKYO-KU
TOKYO 113-0033
JAPAN
E-MAIL: takemura@stat.t.u-tokyo.ac.jp

R. J. ADLER
FACULTY OF INDUSTRIAL ENGINEERING
  AND MANAGEMENT
TECHNION
HAIFA 32000
ISRAEL
E-MAIL: robert@ieadler.technion.ac.il